\documentclass[12pt]{article}

\usepackage{latexsym, amssymb, amsmath, amscd, amsfonts, epsfig, graphicx, colordvi,verbatim,ifpdf}
\usepackage{amsfonts, amsmath, amssymb}
\usepackage{amssymb,amsfonts,amsmath,latexsym,epsfig,cite, psfrag,eepic,color}
\usepackage{amscd,graphics}
\usepackage{latexsym, amssymb,  amsmath,amscd, amsfonts, epsfig, graphicx, colordvi,amsthm}

\usepackage{graphicx}
\usepackage{color}
\usepackage{ifpdf}
\usepackage{fancybox}

\newtheorem{thm}{Theorem}[section]

\newtheorem{defi}[thm]{Definition}

\newtheorem{coro}[thm]{Corollary}

\def\pf{\noindent{\it Proof.} }
\def\qed{\nopagebreak\hfill{\rule{4pt}{7pt}}
\medbreak}

\numberwithin{equation}{section}

\def\qed{\nopagebreak\hfill{\rule{4pt}{7pt}}
\medbreak}

\newlength{\boxedparwidth}
  {\begin{center} \begin{tabular}{|@{\hspace{.315in}}c@{\hspace{.15in}}|}
                  \hline \\ \begin{minipage}[t]{\boxedparwidth}
                  \setlength{\parindent}{.25in}}%
  {\end{minipage} \\ \\ \hline \end{tabular} \end{center}}

\begin{document}
\begin{center}

{ \large\bf A crank for bipartitions with designated summands}
\end{center}

\vskip 5mm

\begin{center}
{  Robert X. J. Hao}$^{1}$, {  Erin Y. Y. Shen}$^{2}$  \vskip 2mm

 $^{1}$ Department of Mathematics and Physics, \\Nanjing Institute of Technology\\
    Nanjing 211167, P.~R.~China\\[6pt]
    $^{2}$School of Science, Hohai University\\
    Nanjing 210098, P.~R.~China\\[6pt]

     \vskip 2mm

    $^1$ haoxj@njit.edu.cn, $^2$shenyy@hhu.edu.cn

\end{center}

\vskip 6mm \noindent {\bf Abstract.}
Andrews, Lewis and Lovejoy introduced the
partition function $PD(n)$ as the number of partitions of $n$ with designated summands. A bipartition of $n$ is an ordered pair of partitions $(\pi_1, \pi_2)$ with the sum of all of the parts being $n$.
In this paper,  we introduce a generalized crank named
the $pd$-crank for bipartitions with designated summands and give some
inequalities for the $pd$-crank of bipartitions with designated summands
modulo 2 and 3. We also define the $pd$-crank moments weighted by
the parity of $pd$-cranks $\mu_{2k,bd}(-1,n)$ and show the positivity
 of $(-1)^n\mu_{2k,bd}(-1,n)$. Let $M_{bd}(m,n)$ denote the number of bipartitions of $n$ with designated summands with $pd$-crank $m$. We prove a monotonicity  property of $pd$-cranks of bipartitions with designated summands and find that the sequence $\{M_{bd}(m,n)\}_{|m|\leq n}$
is unimodal for $n\not= 1,5,7$.

\noindent {\bf Keywords}: bipartition with designated summands; $pd$-crank;
moment; monotonicity

\noindent {\bf AMS Classifications}: 11P81, 11P83, 05A17, 05A20

\section{\large Introduction}

A partition $\lambda$ of a positive integer $n$
is a finite nonincreasing sequence of positive integers $\lambda=(\lambda_1,\ldots, \lambda_l)$ such that
$|\lambda|=\sum_{i=1}^{l}\lambda_i=n$. Let $p(n)$ denote the number of
partitions of $n$. There are two vital statistics in the theory
of partitions named Dyson's rank\cite{Dyson-1944} and the Andrews--Garvan--Dyson crank\cite{Andrews-Garvan-1988}.
The rank \cite{Dyson-1944} of a partition $\lambda$,
denoted by $r(\lambda)$, is defined as the largest part of $\lambda$
minus the number of the parts. The rank can be used to provide combinatorial
interpretations for the following two Ramanujan's famous congruences
\begin{align}
\label{eqnr1}p(5n+4)&\equiv0 \pmod{5}, \\[5pt]
\label{eqnr2}p(7n+5)&\equiv0 \pmod{7}.
\end{align}
The crank\cite{Andrews-Garvan-1988} of a partition $\lambda$ of $n>1$
is defined as
\[c(\lambda)=\left\{
\begin{array}{ll}
\lambda_1,\ \ &\text{ if}\ \ n_1(\lambda)=0,\\[5pt]
\mu(\lambda)-n_1(\lambda), \ \ &\text{ if }\ n_1(\lambda)>0,
\end{array}\right.
\]
where $n_1(\lambda)$ denotes the number of parts equal to
one in $\lambda$ and $\mu(\lambda)$ denotes the number of
parts in $\lambda$ larger than $n_1(\lambda)$.
The crank can be used to provide combinatorial
interpretations for congruences \eqref{eqnr1}, \eqref{eqnr2} as well as
\begin{align}
\label{eqnr3}p(11n+6)&\equiv0\pmod{11}.
\end{align}
Let $M(m,n)$ denote the number of partitions with
crank $m$. For $n\leq 1$ and $m\neq \pm 1, 0$, we set $M(m,n)=0$. For $n=1$
 and $m=\pm 1, 0$, we define
\[M(-1,1)=M(1,1)=1,\  M(0,1)=-1.\]
The generating function of $M(m,n)$ is obtained by Andrews and Garvan\cite{Andrews-Garvan-1988, Garvan-1988} as given by
\begin{align}
C(z;q)=\sum_{m=-\infty}^{\infty} \sum_{n=0}^{\infty} M(m,n)z^mq^n
=\frac{(q;q)_\infty}{(zq;q)_\infty(z^{-1}q;q)_\infty}, \label{gf-c-t-2}
\end{align}
where here and throughout this paper,
$(a;q)_\infty$ stands for the $q$-shifted factorial
\[
(a;q)_{\infty}=\prod_{n=1}^{\infty}(1-aq^{n-1}),\,\, |q|<1.
\]
Let $M(m,t,n)$ denote the number of
partitions of $n$ with crank congruent to $m$ modulo $t$.
Setting $z=-1$ in \eqref{gf-c-t-2}, we have
\begin{align}
C(-1,q)=\sum_{n=0}^{\infty}(M(0,2,n)-M(1,2,n))q^n
=\frac{(q;q)_{\infty}}{(-q;q)^2_{\infty}}.
\end{align}
Andrews and Lewis
\cite{Andrews-Lewis-2000} proved that
\begin{align}
M(0,2,2n)&>M(1,2,2n),\label{inequal1}\\[3pt]
M(1,2,2n+1)&>M(0,2,2n+1),\label{inequal2}
\end{align}
for $n\geq 0$. The above two inequalities imply that the signs of the coefficients of
$\displaystyle \frac{(q;q)_{\infty}}{(-q;q)^2_{\infty}}$ are alternating.

Later, to study the higher-order spt-function $spt_k(n)$, Garvan \cite{Garvan-2011} introduced the $k$-th symmetrized moments $\mu_k(n)$ of cranks of partitions of $n$.  He defined the $k$-th symmetrized moments $\mu_k(n)$ of cranks of partitions of $n$ as
\begin{align}\label{symCrankMom}
\mu_{k}(n)=\sum_{m=-\infty}^{\infty}{m+\lfloor\frac{k-1}{2}\rfloor \choose k}M(m,n).
\end{align}
 Since the fact that $M(m,n)=M(-m,n)$\cite[(1.9)]{Andrews-Garvan-1988}, it is easy to see that $\mu_{2k+1}(n)=0$.
In a recent work, Ji and Zhao\cite{Ji-Zhao-2017} introduced the
crank moments weighted by the parity of cranks as given by
\begin{align}\label{weiCrankMom}
\mu_{2k}(-1,n)=\sum_{m=-\infty}^{\infty}
{m+k-1 \choose 2k}(-1)^mM(m,n).
\end{align}
Let $k=0$ in \eqref{weiCrankMom}, we get
\[
\mu_{0}(-1,n)=M(0,2,n)-M(1,2,n).
\]
Ji and Zhao\cite{Ji-Zhao-2017} proved
the following positivity property of  $(-1)^n\mu_{2k}(-1,n)$.

\begin{thm}{\rm (\!\!\cite[Theorem 1.2]{Ji-Zhao-2017})}
For all $n\geq k\geq 0$, we have
\begin{align}
(-1)^n\mu_{2k}(-1,n)>0.
\end{align}
\end{thm}
This property implies the two inequalities \eqref{inequal1} and \eqref{inequal2} of Andrews and Lewis\cite{Andrews-Lewis-2000}.

Recently, Ji and Zang\cite{Ji-Zang-2018} proved the following
 monotonicity  property of cranks of partitions.

\begin{thm}{\rm (\!\!\cite[Theorem 1.7]{Ji-Zang-2018})}\label{n-1}
For all $n\geq 44$ and $1\leq m \leq n-1$, we have
\begin{align}
M(m-1,n)\geq M(m,n).
\end{align}
\end{thm}
Considering Theorem \ref{n-1} and the symmetry  $M(m,n)=M(-m,n)$,
the following corollary is true.
\begin{coro}{\rm (\!\!\cite[Corollary 1.8]{Ji-Zang-2018})}
For $n\geq 44$, we have
\[
M(1-n,n)\leq\cdots\leq M(-1,n)\leq M(0,n)
\geq M(1,n)\geq\cdots \geq M(n-1,n).
\]
\end{coro}
This means the sequence $\{M(m,n)\}_{|m|\leq n-1}$
is unimodal for $n\geq 44$.

In this paper, we aim to study  bipartitions with designated summands.
In\cite{Andrews-Lewis-Lovejoy-2002}, Andrews, Lewis and Lovejoy
 studied the number of
 partitions with designated summands which are defined on ordinary partitions by tagging exactly one part of each part size. Let $PD(n)$ denote the number of partitions of $n$ with designated summands.  There are  ten partitions of $4$ with designated summands:\emph{}
\begin{align*}
\begin{array}{ccccc}
4',      & 3'+1',  & 2'+2,    & 2+2',    & 2'+1'+1,  \\[5pt]
2'+1+1', &1'+1+1+1,& 1+1'+1+1,& 1+1+1'+1,& 1+1+1+1'.
\end{array}
\end{align*}
Hence we have $PD(4)=10$.
Andrews, Lewis and Lovejoy\cite{Andrews-Lewis-Lovejoy-2002}
derived the generating function of $PD(n)$ as given by
\begin{align}\label{PD(n)}
\sum_{n=0}^{\infty}PD(n)q^n=\frac{(q^6;q^6)_{\infty}}{(q;q)_{\infty}
(q^2;q^2)_{\infty}(q^3;q^3)_{\infty}}.
\end{align}
Andrews, Lewis and Lovejoy\cite{Andrews-Lewis-Lovejoy-2002} obtained
a Ramanujan-type congruence for the partition function $PD(n)$ as
\begin{align}
\label{modPD(3n+2)}PD(3n+2)\equiv 0\pmod{3}.
\end{align}
By introducing the $pd$-rank for partitions with designated summands, Chen, Ji, Jin and the second author\cite{Chen-Ji-Jin-Shen-2013} gave a combinatorial interpretation of the congruence \eqref{modPD(3n+2)}.

A bipartition $\pi$ of $n$ is an ordered pair of partitions
$(\pi_1, \pi_2)$ with  $|\pi_1|+|\pi_2|=n$.
Let $p_{-2}(n)$ denote the number of bipartitions of $n$.
The generating function of $p_{-2}(n)$ is
\[
\sum_{n=0}^{\infty}p_{-2}(n)q^n=\frac{1}{(q;q)^2_\infty}.
\]
A bipartition with designated summands means a bipartition
$\pi=(\pi_1,\pi_2)$ for which $\pi_1$ and $\pi_2$ are both
partitions with designated summands. Here $\pi_1$ and $\pi_2$ are allowed to have one part size tagged in common. Let $PD_{-2}(n)$
denote the number of bipartitions of $n$ with designated summands.
We get the generating function of  $PD_{-2}(n)$ as given by
\begin{align}\label{PDPD(n)}
\sum_{n=0}^{\infty}PD_{-2}(n)q^n=\frac{(q^6;q^6)^2_{\infty}}{(q;q)^2_{\infty}
(q^2;q^2)^2_{\infty}(q^3;q^3)^2_{\infty}}.
\end{align}
Recently, arithmetic properties of  bipartitions with designated summands
have drawn a number of interest, see, for example\cite{Naika-Shivashankar-2016, Hao-shen-2018}.

The main objective of this paper is to investigate bipartitions with designated summands from three aspects. First, we introduce a generalized crank named
the $pd$-crank for bipartitions with designated summands and establish some
inequalities for the $pd$-crank of bipartitions with designated summands
modulo 2 and 3.  Second, we define the $pd$-crank moments weighted by
the parity of $pd$-cranks $\mu_{2k,bd}(-1,n)$ and show the positivity
property of $(-1)^n\mu_{2k,bd}(-1,n)$. Finally, we prove a monotonicity  property of $pd$-cranks of bipartitions with designated summands.

\section{\large The $pd$-crank and its inequalities}

In this section, we first introduce a generalized
crank which called $pd$-crank for bipartitions with designated summands.
The definition of the $pd$-crank relies on the construction of the following bijection $\Delta$ which Chen, Ji, Jin and the second
author established in \cite{Chen-Ji-Jin-Shen-2013}.

\begin{thm}\label{lambdaab}{\rm(\!\!\cite[Theorem 3.1]{Chen-Ji-Jin-Shen-2013})} There is a
bijection $\Delta$ between the set of partitions of $n$ with designated summands  and the set of vector partitions  $(\alpha, \beta)$  with $|\alpha|+|\beta|=n$,
where $\alpha$ is an ordinary partition and $\beta$  is a partition into parts $\not\equiv
\pm1 \pmod{6}$.
\end{thm}

Chen, Ji, Jin and the second
author gave a combinatorial proof of the above theorem
 which illustrates the construction of the bijection $\Delta$, see \cite[Theorem 3.1]{Chen-Ji-Jin-Shen-2013}.
Under the map  $\Delta$, Chen, Ji, Jin and the second author defined
 the $pd$-rank of a partition $\lambda$ with designated summands
 in terms of the pair of partitions  $(\alpha, \beta)$.

\begin{defi}{\rm(\!\!\cite[Definition 3.2]{Chen-Ji-Jin-Shen-2013})}
Let $\lambda$ be a partition with designated summands and let $(\alpha,\beta)=\Delta(\lambda)$.   Then the $pd$-rank of  $\lambda$,  denoted  $r_d(\lambda)$, is defined by
\begin{align}
r_d(\lambda)=l_e(\alpha)-l_e(\beta),
\end{align}
where
$l_e(\alpha)$ is the number of even parts of $\alpha$ and
and $l_e(\beta)$ is the number of even parts of $\beta$.
\end{defi}
We are now ready to give the definition of the $pd$-crank for
a bipartition $(\lambda_1,\lambda_2)$ with designated summands under the map
$\Delta$.
\begin{defi}\label{defpd}
Let $(\lambda_1,\lambda_2)$ be a bipartition with
designated summands and let
$$
\left((\alpha_1,\beta_1),(\alpha_2,\beta_2))
=(\Delta(\lambda_1),\Delta(\lambda_2)\right).
$$
Then the $pd$-crank of  $(\lambda_1,\lambda_2)$,  denoted  $r_{bd}(\lambda)$, is defined by
\begin{align}
r_{bd}(\lambda)=l(\alpha_1)-l(\alpha_2),
\end{align}
where
$l(\alpha)$ is the number of parts of $\alpha$.
\end{defi}

Let $M_{bd}(m,n)$ denote the number of bipartitions of $n$ with designated summands with $pd$-crank $m$,
and $M_{bd}(k,t,n)$ denote the number of bipartitions with designated summands of $n$ with $pd$-crank congruent to $k$ modulo $t$.
According to the Definition \ref{defpd}, it is clear that
\begin{align}\label{pdsym}
M_{bd}(m,n)=M_{bd}(-m,n).
\end{align}
Following the results of Andrews and Lewis\cite{Andrews-Lewis-2000}, we
investigate the inequalities of the $pd$-crank of bipartitions with designated summands modulo $2$ and $3$. By the definition of the $pd$-crank of
bipartitions with designated summands, we have
\begin{align}\label{Nbdz3}
\sum_{m=-\infty}^{\infty}\sum_{n=0}^{\infty}M_{bd}(m,n)z^mq^n
=\frac{(q^6;q^6)^2_{\infty}}{(zq;q)_{\infty}(z^{-1}q;q)_{\infty}
(q^2;q^2)^2_{\infty}(q^3;q^3)^2_{\infty}}.
\end{align}
The following theorem shows that the sequence $M_{bd}(0,2,n)-M_{bd}(1,2,n)$
alternates in sign.
\begin{thm}\label{thmmod2}
For $n\geq 0$, we have
\begin{align}
M_{bd}(0,2,2n)&>M_{bd}(1,2,2n),\\
M_{bd}(1,2,2n+1)&>M_{bd}(0,2,2n+1).
\end{align}
\end{thm}

\pf
Setting $z=-1$ in \eqref{Nbdz3}, we get
\begin{align}\label{N02}
\sum_{n=0}^{\infty}(M_{bd}(0,2,n)-M_{bd}(1,2,n))q^n
&=\frac{(q^6;q^6)^2_{\infty}}{(-q;q)_{\infty}^2
(q^2;q^2)^2_{\infty}(q^3;q^3)^2_{\infty}}\nonumber\\[3pt]
&=\frac{(q^6;q^6)^2_{\infty}(q;q)^2_{\infty}}{
(q^2;q^2)^4_{\infty}(q^3;q^3)^2_{\infty}}=\frac{f^2_6f^2_1}{f^4_2f^2_3},
\end{align}
where here and throughout this paper,
$f_k$ is defined by
\[
f_k=(q^k;q^k)_{\infty}.
\]
By the following $2$-dissection{\rm(\!\!\cite[Lemma 2.6]{Yao-Xia-2013})}
\begin{align*}
\frac{f_1}{f_3}&=\frac{f_2f_{16}f^2_{24}}
{f^2_6f_8f_{48}}-q\frac{f_2f^2_8f_{12}f_{48}}
{f_4f^2_6f_{16}f_{24}},
\end{align*}
we can deduce that
\begin{align}\label{2f2f3}
\frac{f_1^2}{f_3^2}&=\frac{f_2^2f_{16}^2f^4_{24}}
{f^4_6f_{8}^2f_{48}^2}+q^2\frac{f_2^2f^4_8f_{12}^2f_{48}^2}
{f_4^2f^4_6f_{16}^2f_{24}^2}-
2q\frac{f_2^2f_8f_{12}f_{24}}
{f_4f^4_6}.
\end{align}
Applying \eqref{2f2f3} into \eqref{N02}, we get
\begin{align*}
\sum_{n=0}^{\infty}(M_{bd}(0,2,n)-M_{bd}(1,2,n))q^n
&=\frac{f_6^2}{
f^4_2}\left(\frac{f_2^2f_{16}^2f^4_{24}}
{f^4_6f_{8}^2f_{48}^2}+q^2\frac{f_2^2f^4_8f_{12}^2f_{48}^2}
{f_4^2f^4_6f_{16}^2f_{24}^2}-
2q\frac{f_2^2f_8f_{12}f_{24}}
{f_4f^4_6}\right),
\end{align*}
which implies the coefficient of $q^n$ in \eqref{N02} is positive(negative) when $n$ is even(odd). This completes the proof.
\qed

Considering the $pd$-crank of bipartitions with designated summands modulo $3$, we have the following results.

\begin{thm}\label{thmmod3}
For $n\geq 0$, we have
\begin{align}
M_{bd}(0,3,6n)&>M_{bd}(1,3,6n),\label{6n} \\
M_{bd}(0,3,6n+1)&<M_{bd}(1,3,6n+1),\label{6n+1} \\
M_{bd}(0,3,6n+2)&>M_{bd}(1,3,6n+2),\label{6n+2} \\
M_{bd}(0,3,6n+3)&>M_{bd}(1,3,6n+3),\label{6n+3} \\
M_{bd}(0,3,6n+4)=M_{bd}(1,&3,6n+4)=M_{bd}(2,3,6n+4),\label{6n+4}\\
M_{bd}(0,3,6n+5)&<M_{bd}(1,3,6n+5).\label{6n+5}
\end{align}
\end{thm}

\pf
Substituting $z=\zeta=e^{\frac{2\pi i}{3}}$ into \eqref{Nbdz3},
we find that
\begin{align}
\sum_{m=-\infty}^{\infty}\sum_{n=0}^{\infty}M_{bd}(m,n)\zeta^mq^n
&=\sum_{n=0}^{\infty}\sum_{k=0}^{2}M_{bd}(k,3,n)\zeta^kq^n
\nonumber \\[5pt]
&=\sum_{n=0}^{\infty}(M_{bd}(0,3,n)-M_{bd}(1,3,n))q^n\nonumber \\[3pt]
&=\frac{(q^6;q^6)^2_{\infty}}
{(\zeta q;q)_{\infty}(\zeta^{-1}q;q)_{\infty}
(q^2;q^2)^2_{\infty}(q^3;q^3)^2_{\infty}}
\nonumber \\[5pt]
\label{pro1}&=\frac{(q;q)_{\infty}(q^6;q^6)^2_{\infty}}
{(q^2;q^2)^2_{\infty}(q^3;q^3)^3_{\infty}}.
\end{align}
Using the $2$-dissection of $\displaystyle \frac{(q;q)_\infty}{(q^3;q^3)^3_{\infty}}$
\cite[Lemma 2.5]{Yao-Xia-2013}, which is
\begin{align}
\frac{(q;q)_\infty}{(q^3;q^3)^3_{\infty}}
&=\frac{(q^2;q^2)_{\infty}(q^4;q^4)^2_{\infty}
(q^{12};q^{12})^2_{\infty}}
{(q^6;q^6)^7_{\infty}}-q\frac{(q^2;q^2)^3_{\infty}
(q^{12};q^{12})^6_{\infty}}{(q^4;q^4)^2_{\infty}(q^6;q^6)^9_{\infty}},
\end{align}
we get
\begin{align}\label{dis2}
&\sum_{n=0}^{\infty}(M_{bd}(0,3,n)-M_{bd}(1,3,n))q^n\nonumber\\
&\hskip 1cm =\frac{(q^6;q^6)^2_{\infty}}{(q^2;q^2)^2_{\infty}}
\left(\frac{(q^2;q^2)_{\infty}(q^4;q^4)^2_{\infty}
(q^{12};q^{12})^2_{\infty}}
{(q^6;q^6)^7_{\infty}}-q\frac{(q^2;q^2)^3_{\infty}
(q^{12};q^{12})^6_{\infty}}{(q^4;q^4)^2_{\infty}(q^6;q^6)^9_{\infty}}\right).
\end{align}
Therefore we have
\begin{align}\label{even}
\sum_{n=0}^{\infty}(M_{bd}(0,3,2n)-M_{bd}(1,3,2n))q^n
=\frac{(q^2;q^2)^2_{\infty}(q^6;q^6)^2_{\infty}}
{(q;q)_\infty(q^3;q^3)^5_{\infty}}
=\frac{(q^6;q^6)^2_{\infty}}{(q^3;q^3)^5_{\infty}}\psi(q).
\end{align}
By the identity \cite[p. 49]{Berndt-1991}
\[\psi(q)=f(q^3,q^6)+q\psi(q^9),\]
we deduce that
\begin{align}\label{pro3}
\sum_{n=0}^{\infty}(M_{bd}(0,3,2n)-M_{bd}(1,3,2n))q^n
=\frac{(q^6;q^6)^2_{\infty}}{(q^3;q^3)^5_{\infty}}
\left(f(q^3,q^6)+q\psi(q^9)\right).
\end{align}
Since the coefficients of $q^{3n}$ and $q^{3n+1}$ in \eqref{pro3} are both positive, \eqref{6n} and \eqref{6n+2} are true. Noting that the coefficient of $q^{3n+2}$ in \eqref{pro3} is zero, hence we have \eqref{6n+4}.

According to \eqref{dis2},  we find that
\begin{align}\label{odd1}
\sum_{n=0}^{\infty}(M_{bd}(1,3,2n+1)-M_{bd}(0,3,2n+1))q^n
&=\frac{(q;q)_{\infty}
(q^{6};q^{6})^6_{\infty}}{(q^2;q^2)^2_{\infty}(q^3;q^3)^7_{\infty}}\nonumber\\
&=\frac{1}{\psi(q)}\cdot\frac{(q^{6};q^{6})^6_{\infty}}
{(q^3;q^3)^7_{\infty}}.
\end{align}
Let
\begin{align*}
A(q)=\frac{f_2f^2_3}{f_1f_6},
\end{align*}
by the identity \cite[Lem. 2.2]{Hao-shen-2018}
\begin{align}
\label{psiq3}\frac{1}{\psi(q)}=\frac{\psi(q^9)}{\psi(q^3)^4}
\left(A(q^3)^2-qA(q^3)\psi(q^9)+q^2\psi(q^9)^2\right),
\end{align}
we may obtain that
\begin{align}\label{odd2}
&\sum_{n=0}^{\infty}(M_{bd}(1,3,2n+1)-M_{bd}(0,3,2n+1))q^n\nonumber\\
&\hskip 1.5cm =\frac{(q^{6};q^{6})^6_{\infty}}
{(q^3;q^3)^7_{\infty}}\cdot\frac{\psi(q^9)}{\psi(q^3)^4}
\left(A(q^3)^2-qA(q^3)\psi(q^9)+q^2\psi(q^9)^2\right).
\end{align}
 It is clear that the coefficients of $q^{3n}$ and $q^{3n+2}$
 in \eqref{odd2} are both positive and the coefficient of $q^{3n+1}$
 is negative. These lead to \eqref{6n+1}, \eqref{6n+5} and \eqref{6n+3}
 respectively. This completes the proof.\qed

Theorem \ref{thmmod3} shows that the $pd$-crank  we defined can be
used to divide the set of bipartitions of $6n+4$ with designated summands
into three equinumerous classes. Hence we provide a combinatorial
interpretation for the congruence
\begin{align}
PD_{-2}(6n+4)\equiv 0 \pmod{3}
\end{align}
proved by Mahadeva Naika and Shivashankar \cite[(3.4)]{Naika-Shivashankar-2016}.

\section{\large The $pd$-crank moments weighted by the parity of $pd$-cranks}

As a natural analog to the $2k$-th crank
moments $\mu_{2k}(-1,n)$ weighted by the parity of cranks due to Ji and Zhao\cite{Ji-Zhao-2017}, we define the $2k$-th $pd$-crank moments of bipartitions with designated summands weighted by the
parity of $pd$-cranks as given by
\begin{align}
\mu_{2k,bd}(-1,n)=\sum_{m=-\infty}^{\infty}{m+k-1 \choose 2k}
(-1)^mM_{bd}(m,n).
\end{align}
In this section, we aim to prove the following positivity property of $(-1)^n\mu_{2k,bd}(-1,n)$.

\begin{thm}\label{mubd}
For $n\geq k\geq 0$, we have
\begin{align}\label{muinequ}
(-1)^n\mu_{2k,bd}(-1,n)>0.
\end{align}
\end{thm}

It is worth mentioning that Theorem \ref{mubd} reduces to
Theorem \ref{thmmod2} when $k=0$ in \eqref{muinequ}.

With the aid of the proof of Theorem 1.3 of Ji and Zhao\cite{Ji-Zhao-2017} and applying Andrews' $j$-fold generalization of $q$-Whipple's theorem
\cite{Andrews-1975}, we obtain the generating function of $\mu_{2k,bd}(-1,n)$
as
\begin{align}\label{gfun-1}
&\sum_{n=0}^{\infty}\mu_{2k,bd}(-1,n)q^n \nonumber\\
&\hskip 0.5cm=\frac{(q^6;q^6)^2_{\infty}}{(q^2;q^2)^2_{\infty}(q^3;q^3)^2_{\infty}
(-q;q)^2_{\infty}}\sum_{ n_{k}\geq n_{k-1} \cdots \geq n_1\geq
1}\frac{(-1)^k q^{n_1+n_2+\cdots+{n_{k}}}}
{(1+q^{n_1})^2(1+q^{n_2})^2\cdots(1+q^{n_{k}})^2}.
\end{align}
Considering the proof of Theorem 1.4 of Ji and Zhao\cite{Ji-Zhao-2017},
we find that the generating function  \eqref{gfun-1} is equivalent to
\begin{align}\label{gfun-2}
&\sum_{n=0}^{\infty}\mu_{2k,bd}(-1,n)q^n \nonumber\\
&\hskip 0.5cm=\frac{(q^6;q^6)^2_{\infty}}{(q^2;q^2)^2_{\infty}(q^3;q^3)^2_{\infty}
(-q;q)^2_{\infty}}\sum_{ m_{k}> m_{k-1}> \cdots >m_1\geq
1}\frac{(-1)^{m_k}m_1(m_2-m_1)\cdots(m_k-m_{k-1}) q^{m_k}}
{(1-q^{m_1})(1-q^{m_2})\cdots(1-q^{m_{k}})}.
\end{align}

We are now stand at the point to prove Theorem \ref{mubd}.

\emph{Proof of Theorem\ref{mubd}.}
Replacing $q$ by $-q$ in \eqref{gfun-2}, we get
\begin{align}\label{gfun-3}
&\sum_{n=0}^{\infty}(-1)^n\mu_{2k,bd}(-1,n)q^n \nonumber\\
&\hskip 0.5cm=\frac{(q^6;q^6)^2_{\infty}}{(q^2;q^2)^2_{\infty}(-q^3;-q^3)^2_{\infty}
(q;-q)^2_{\infty}}\sum_{ m_{k}> m_{k-1}> \cdots >m_1\geq
1}\frac{m_1(m_2-m_1)\cdots(m_k-m_{k-1}) q^{m_k}}
{(1-(-q)^{m_1})(1-(-q)^{m_2})\cdots(1-(-q)^{m_{k}})}.
\end{align}
Using the facts that
\begin{align*}
\frac{1}{(q;-q)_{\infty}}&=(-q;q^2)_{\infty},\\[5pt]
(-q^3;-q^3)_{\infty}&=(-q^3;q^6)_{\infty}(q^6;q^6)_{\infty},
\end{align*}
\eqref{gfun-3} turns out to be
\begin{align}\label{gfun-4}
&\sum_{n=0}^{\infty}(-1)^n\mu_{2k,bd}(-1,n)q^n \nonumber\\
&\hskip 0.5cm=\frac{(q^6;q^6)_{\infty}(-q;q^2)_{\infty}}
{(q^2;q^2)^2_{\infty}(q^6;q^6)_{\infty}(-q^3;q^6)^2_{\infty}}
\sum_{ m_{k}> m_{k-1}> \cdots >m_1\geq
1}\frac{(-q;q^2)_{\infty}\cdot m_1(m_2-m_1)\cdots(m_k-m_{k-1}) q^{m_k}}
{(1-(-q)^{m_1})(1-(-q)^{m_2})\cdots(1-(-q)^{m_{k}})}\nonumber \\
&\hskip 0.5cm=\frac{(q^6;q^6)_{\infty}(-q;q^6)_{\infty}(-q^5;q^6)_{\infty}}
{(q^2;q^2)^2_{\infty}(q^6;q^6)_{\infty}(-q^3;q^6)_{\infty}}
\sum_{ m_{k}> m_{k-1}> \cdots >m_1\geq
1}\frac{(-q;q^2)_{\infty}\cdot m_1(m_2-m_1)\cdots(m_k-m_{k-1}) q^{m_k}}
{(1-(-q)^{m_1})(1-(-q)^{m_2})\cdots(1-(-q)^{m_{k}})}
\end{align}
Multiplying the right hand side of \eqref{gfun-4}  by
\[\frac{(q^3;q^6)_{\infty}(-q^2;q^6)_{\infty}(-q^4;q^6)_{\infty}}
{(q^3;q^6)_{\infty}(-q^2;q^6)_{\infty}(-q^4;q^6)_{\infty}},\]
we get
\begin{align}\label{gfun-6}
&\sum_{n=0}^{\infty}(-1)^n\mu_{2k,bd}(-1,n)q^n \nonumber\\
&\hskip 0.5cm=\frac{(-q;q^3)_{\infty}(-q^2;q^3)_{\infty}(q^3;q^3)_{\infty}}
{(q^2;q^2)_{\infty}(q^6;q^6)^2_{\infty}(q^4;q^{12})_{\infty}
(q^6;q^{12})_{\infty}(q^8;q^{12})_{\infty}}
\nonumber \\[5pt]
&\hskip 1.5cm \times\sum_{ m_{k}> m_{k-1}> \cdots >m_1\geq
1}\frac{(-q;q^2)_{\infty}\cdot m_1(m_2-m_1)\cdots(m_k-m_{k-1}) q^{m_k}}
{(1-(-q)^{m_1})(1-(-q)^{m_2})\cdots(1-(-q)^{m_{k}})}.
\end{align}
Applying the Jacobi triple product identity\cite[(1.3.10)]{Berndt-2006}
\begin{align*}
(-z;q)_\infty({-q/z};q)_\infty(q;q)_\infty
=\sum_{n=-\infty}^{\infty}z^{n}q^{n\choose 2}
\end{align*}
into \eqref{gfun-6} with $q$ replaced by $q^3$
and $z$ replaced by $q$, we deduce that
\begin{align}\label{gfun-5}
&\sum_{n=0}^{\infty}(-1)^n\mu_{2k,bd}(-1,n)q^n \nonumber\\
&\hskip 0.5cm=\frac{\sum_{n=-\infty}^{\infty}q^{n(3n-1)/2}}
{(q^2;q^2)_{\infty}(q^6;q^6)^2_{\infty}(q^4;q^{12})_{\infty}
(q^6;q^{12})_{\infty}(q^8;q^{12})_{\infty}
}\nonumber\\[5pt]
&\hskip 1.5cm \times
\sum_{ m_{k}> m_{k-1}> \cdots >m_1\geq
1}\frac{(-q;q^2)_{\infty}\cdot m_1(m_2-m_1)\cdots(m_k-m_{k-1}) q^{m_k}}
{(1-(-q)^{m_1})(1-(-q)^{m_2})\cdots(1-(-q)^{m_{k}})}
\end{align}
Given $m_{k}> m_{k-1}> \cdots >m_1\geq 1$, we define
\begin{align}\label{sum}
\sum_{m\geq 0}h_{m_1,m_2,\ldots,m_k}(m)q^m
=\frac{(-q;q^2)_{\infty}}
{(1-(-q)^{m_1})(1-(-q)^{m_2})\cdots(1-(-q)^{m_{k}})}.
\end{align}
In \eqref{sum}, when $m_i$ is odd, the corresponding term
$(1-(-q)^{m_i})$ in denominator can be canceled by
$(-q;q^2)_{\infty}$  since $m_i$ differs from each other.
When $m_i$ is even, the corresponding term
$1/(1-(-q)^{m_i})$ does not have negative coefficients.
Therefore we find that $h_{m_1,m_2,\ldots,m_k}(m)\geq 0$
and $h_{m_1,m_2,\ldots,m_k}(0)=1$. This implies that
$(-1)^n\mu_{2k,bd}(-1,n)>0$ for any nonnegative integer $n\geq k$.
This completes the proof.\qed

\section{\large A monotonicity  property of $M_{bd}(m,n)$}

In this section, we investigate the following  monotonicity property of $M_{bd}(m,n)$ which leads to the unimodality of the sequence $\{M_{bd}(m,n)\}_{|m|\leq n}$ with $n\not= 1,5,7$.

\begin{thm}\label{thmuni}
For  $m\geq 1$ and $n\geq 0$, \begin{equation}
M_{bd}(m-1,n)\geq M_{bd}(m,n),
\end{equation}
except for $(m,n)=(1,1),(1,5)$ or $(1,7)$.
\end{thm}

Fu and Tang \cite{Fu-Tang-2018} studied a generalized crank named $k$-crank for $k$-colored partitions. Let $M_k(m,n)$ denote the number of $k$-colored partitions of $n$ with $k$-crank $m$.  The generating function of $M_k(m,n)$ is given by
\begin{equation}\label{M2}
\sum_{m=-\infty}^{\infty}\sum_{n=0}^{\infty}
M_k(m,n)z^mq^n=\frac{(q;q)_{\infty}^{2-k}}
{(zq;q)_{\infty}(z^{-1}q;q)_{\infty}}
\end{equation}
for $k\in N$.

Now we provide a proof for Theorem \ref{thmuni}.

\emph{Proof of Theorem \ref{thmuni}.}
Setting $k=2$ in \eqref{M2} and applying it to \eqref{Nbdz3},
we get
\begin{align*}
\sum_{m=-\infty}^{\infty}\sum_{n=0}^\infty M_{bd}(m,n)z^mq^n&=
\frac{(q^6;q^6)_{\infty}^2}{(zq;q)_{\infty}(z^{-1}q;q)_{\infty}
(q^2;q^2)_{\infty}^2(q^3;q^3)_{\infty}^2}\\
&=\frac{(q^6;q^6)_{\infty}^2}{(q^2;q^2)_{\infty}^2(q^3;q^3)_{\infty}^2}
\sum_{m=-\infty}^{\infty}\sum_{n=0}^{\infty}M_2(m,n)z^mq^n.
\end{align*}
It is clear that
\begin{align}\label{Mk}
&\sum_{n=0}^\infty (M_{bd}(m-1,n)-M_{bd}(m,n))q^n\nonumber \\
&\hskip 1.5cm =\frac{(q^6;q^6)_{\infty}^2}{(q^2;q^2)_{\infty}^2(q^3;q^3)_{\infty}^2}
\sum_{n=0}^{\infty}(M_2(m-1,n)-M_2(m,n))q^n.
\end{align}
According to \cite[Theorem 1.4]{Zang-Zhang-2018} and \cite[(6.3)]{Zang-Zhang-2018},  Zang and Zhang proved that
$M_k(m-1, n)\geq M_k(m,n)$ for  $m\geq 2$, $k\geq 2$ and $n\geq 0$.
Hence by \eqref{Mk}, we find that $M_{bd}(m-1, n)\geq M_{bd}(m,n)$
when $m\geq 2$ and $n\geq 0$.

For $m=1$, we have
\begin{align}\label{m=1}
&\sum_{n=0}^\infty (M_{bd}(0,n)-M_{bd}(1,n))q^n\nonumber \\
&\hskip 1.5cm =\frac{(q^6;q^6)_{\infty}^2}{(q^2;q^2)_{\infty}^2(q^3;q^3)_{\infty}^2}
\sum_{n=0}^{\infty}(M_2(0,n)-M_2(1,n))q^n.
\end{align}
Let $\{c_n\}_{n=0}^{\infty}$ be a sequence of nonnegative integers. By \cite[(6.5)]{Zang-Zhang-2018} and a simple calculation, we can derive that
\begin{equation}\label{M201}
\sum_{n=0}^{\infty}(M_2(0,n)-M_2(1,n))q^n=1-q+q^4+\sum_{n=6}^{\infty}c_nq^n.
\end{equation}
Noting that
\begin{align}\label{f236}
\frac{(q^6;q^6)_{\infty}}{(q^2;q^2)_{\infty}(q^3;q^3)_{\infty}}
&=\frac{(q^6;q^6)_{\infty}(-q^3;q^3)_{\infty}}
{(q^2;q^2)_{\infty}(q^3;q^3)_{\infty}(-q^3;q^3)_{\infty}}
\nonumber \\[3pt]
&=\prod_{i=1}^{\infty}\frac{1+q^{3i}}{1-q^{2i}}
=\prod_{i=1}^{\infty}\left(1+\frac{q^{2i}}{1-q^i}\right),
\end{align}
and substituting \eqref{M201}, \eqref{f236} into \eqref{m=1},
we have
\begin{align}\label{Nbd01}
\sum_{n=0}^\infty (M_{bd}(0,n)-M_{bd}(1,n))q^n
=\left(1-q+q^4+\sum_{n=6}^{\infty}c_nq^n \right)
\prod_{i=1}^{\infty}\left(1+\frac{q^{2i}}{1-q^i}\right)^2.
\end{align}
Next, we aim to show the coefficients of
\begin{align}\label{mainpo}
(1-q+q^4)\prod_{i=1}^{\infty}\left(1+\frac{q^{2i}}{1-q^i}\right)^2
\end{align}
are nonnegative except for $n=1,5,7$.

Let $P_2(n)$ denote the number of partitions of $n$ with each part exists at least twice. Let $P_{24}(n)$ denote the number of partitions of $n$ counted by $P_2(n)$ with  $1$ occurs zero, two or four times. We deduce that
\begin{align}\label{mainf1}
&(1-q+q^4)\prod_{i=1}^{\infty}\left(1+\frac{q^{2i}}{1-q^i}\right)\nonumber\\
&\hskip 1.5cm=\prod_{i=1}^{\infty}\left(1+\frac{q^{2i}}{1-q^i}\right)
-q(1-q^3)\left(1+\frac{q^{2}}{1-q}\right)
\prod_{i=2}^{\infty}\left(1+\frac{q^{2i}}{1-q^i}\right)
\nonumber\\
&\hskip 1.5cm=\prod_{i=1}^{\infty}\left(1+\frac{q^{2i}}{1-q^i}\right)
-q(1+q^2+q^4)
\prod_{i=2}^{\infty}\left(1+\frac{q^{2i}}{1-q^i}\right)
\nonumber\\
&\hskip 1.5cm=1+\sum_{n=1}^{\infty}\left(P_2(n)-P_{24}(n-1)\right)q^n.
\end{align}
After a simple calculation, we find that
\begin{align}\label{mainf2}
&(1-q+q^4)\prod_{i=1}^{\infty}\left(1+\frac{q^{2i}}{1-q^i}\right)\nonumber\\
&\hskip 1.5cm=1-q+q^2+2q^4-q^5+4q^6-q^7+6q^8+8q^{10}+15q^{12}
+\sum_{n=14}^{\infty}a_nq^n.
\end{align}
In light of \eqref{mainf1} and \eqref{mainf2}, we wish to construct an injection $\Omega$ from the set of partitions $\alpha$ counted by $P_{24}(n-1)$ to the set of partitions  $\beta$  enumerated by $P_2(n)$
to prove the coefficients $a_n$ of \eqref{mainf2} are nonnegative for $n\geq 14$.

If a partition $\alpha$ counted by $P_{24}(n-1)$ has 1 as its part,
1 may occur twice or four times. By adding a 1 to $\alpha$ as a part, we get a partition $\beta=\Omega(\alpha)$ enumerated by $P_2(n)$.
Hence $\beta$ has three or five 1s. For instance, let $\alpha=(4,4,3,3,1,1)$ and it belongs to the set of partitions counted by $P_{24}(16)$.
 Under the map $\Omega$,  the corresponding partition $\beta=\Omega(\alpha)$ should be $(4,4,3,3,1,1,1)$ and it belongs to the set of partitions counted by $P_{2}(17)$.

Suppose $\alpha$ is a partition enumerated by $P_{24}(n-1)$ with no
part equals to 1. Let the largest part size of the partition
$\alpha$ be $m$ and the second largest part size be $\alpha_1$.
We have the following cases.
\begin{itemize}
\item
If $m \geq 5$ or $m=3$ and $m$ appears more than twice in $\alpha$, we obtain $\beta$  by rewriting one $m$ in $\alpha$ plus one more $1$ as $m+1$ 1s.
For example, let $$
\alpha=(\underbrace{m,m,\ldots,m,}_{l,\ l>2}\alpha_1,\alpha_1,\ldots,\alpha_n,\alpha_n),
$$
the corresponding partition $\beta=\Omega(\alpha)$ should be
$$
\beta=(\underbrace{m,\ldots,m,}_{l-1}\alpha_1,\alpha_1,\ldots,
\alpha_n,\alpha_n,\underbrace{1,1,\ldots,1}_{m+1}).
$$

\item
If $m \geq 5$  or $m=3$ and $m$ appears exactly  twice in $\alpha$, we obtain $\beta$  by rewriting  two $m$ in $\alpha$ plus one more $1$ as $2m+1$ $1$s. For example, let $$
\alpha=(m,m,\alpha_1,\alpha_1,\ldots,\alpha_n,\alpha_n),
$$
the corresponding partition $\beta=\Omega(\alpha)$ should be
$$
\beta=(\alpha_1,\alpha_1,\ldots,
\alpha_n,\alpha_n,\underbrace{1,1,\ldots,1}_{2m+1}).
$$

\item
If $m=4$ and $4$ appears more than three times or exactly twice in $\alpha$, we obtain $\beta$  by rewriting two $4$s in $\alpha$ plus one more $1$ as nine $1$s. For example, let $$
\alpha=(\underbrace{4,4,\ldots,4,}_{l,\ l=2 {\ \rm or\ } l>3}\alpha_1,\alpha_1,\ldots,\alpha_n,\alpha_n),
$$
the corresponding partition $\beta=\Omega(\alpha)$ should be
$$
\beta=(\underbrace{4,4,\ldots,4,}_{l-2}\alpha_1,\alpha_1,\ldots,
\alpha_n,\alpha_n,\underbrace{1,1,\ldots,1}_{9}).
$$

\item
If $m=4$ and $4$ appears exactly three times in $\alpha$, we obtain $\beta$  by rewriting the three $4$s in $\alpha$ plus one more $1$ as $3,3,3,2,2$. For example, let $$
\alpha=(4,4,4,\alpha_1,\ldots,\alpha_1,\alpha_2,\alpha_2,\ldots),
$$
the corresponding partition $\beta=\Omega(\alpha)$ should be
$$
\beta=(3,3,3,\alpha_1,\ldots,\alpha_1,2,2,\alpha_2,\alpha_2,\ldots).
$$

\item
If $m=2$ and $2$ appears more than or equal to seven times in $\alpha$, we obtain $\beta$  by rewriting five $2$s in $\alpha$ plus one more $1$ as $3,3,3,1,1$. For example, let $$
\alpha=(\underbrace{2,2,\ldots,2,}_{l,\ l\geq 7}),
$$
the corresponding partition $\beta=\Omega(\alpha)$ should be
$$
\beta=(3,3,3,\underbrace{2,2,\ldots,2,}_{l-5}1,1).
$$

\end{itemize}
Therefore the map $\Omega$ is an injection. Hence the coefficients $a_n$ of \eqref{mainf2} are nonnegative for $n\geq 14$. Applying \eqref{mainf2} into \eqref{mainpo}, we have
\begin{align}\label{mainf3}
&(1-q+q^4)\prod_{i=1}^{\infty}\left(1+\frac{q^{2i}}{1-q^i}\right)^2\nonumber\\
&\hskip 1cm=\left(1-q+q^2+2q^4-q^5+4q^6-q^7+6q^8+8q^{10}+15q^{12}
+\sum_{n=14}^{\infty}a_nq^n\right) \prod_{i=1}^{\infty}\left(1+\frac{q^{2i}}{1-q^i}\right)\nonumber\\
&\hskip 1cm=(1+q^4+q^6)(1-q+q^4)
\prod_{i=1}^{\infty}\left(1+\frac{q^{2i}}{1-q^i}\right)\nonumber\\
&\hskip 2cm+\left(q^2+3q^6+5q^8+7q^{10}+15q^{12}
+\sum_{n=14}^{\infty}a_nq^n\right)\prod_{i=1}^{\infty}\left(1+\frac{q^{2i}}{1-q^i}\right).
\end{align}
Substituting \eqref{mainf2} into \eqref{mainf3}, we deduce that
\begin{align}\label{mainf4}
&(1-q+q^4)\prod_{i=1}^{\infty}\left(1+\frac{q^{2i}}{1-q^i}\right)^2\nonumber\\
&\hskip 1cm=(1+q^4+q^6)\left(1-q+q^2+2q^4-q^5+4q^6-q^7+6q^8+8q^{10}+15q^{12}
+\sum_{n=14}^{\infty}a_nq^n\right) \nonumber\\
&\hskip 1.5cm+\left(q^2+3q^6+5q^8+7q^{10}+15q^{12}
+\sum_{n=14}^{\infty}a_nq^n\right)\prod_{i=1}^{\infty}
\left(1+\frac{q^{2i}}{1-q^i}\right)\nonumber\\[5pt]
&\hskip 1cm=(1+q^4+q^6)(-q-q^5-q^7)+(1+q^4+q^6)
(1+q^2+2q^4+4q^6+6q^8\nonumber\\[5pt]
&\hskip 1.5cm
+8q^{10}+15q^{12}
+\sum_{n=14}^{\infty}a_nq^n)+q^2\prod_{i=1}^{\infty}
\left(1+\frac{q^{2i}}{1-q^i}\right)\nonumber\\[5pt]
&\hskip 1.5cm+\left(3q^6+5q^8+7q^{10}+15q^{12}
+\sum_{n=14}^{\infty}a_nq^n\right)\prod_{i=1}^{\infty}
\left(1+\frac{q^{2i}}{1-q^i}\right)\nonumber\\[5pt]
&\hskip 1cm=G_1(q)+G_2(q)+G_3(q)+G_4(q),
\end{align}
where
\begin{align}
\label{G1}G_1(q)&=(1+q^4+q^6)(-q-q^5-q^7)=-q-2q^5-2q^7-q^9-2q^{11}-q^{13},\\
G_2(q)&=(1+q^4+q^6)(1+q^2+2q^4+4q^6+6q^8+8q^{10}+15q^{12}
+\sum_{n=14}^{\infty}a_nq^n),\\
G_3(q)&=q^2\prod_{i=1}^{\infty}
\left(1+\frac{q^{2i}}{1-q^i}\right),\\
G_4(q)&=\left(3q^6+5q^8+7q^{10}+15q^{12}
+\sum_{n=14}^{\infty}a_nq^n\right)\prod_{i=1}^{\infty}
\left(1+\frac{q^{2i}}{1-q^i}\right).
\end{align}
It is clear that the coefficients of the functions $G_2(q)$, $G_3(q)$ and $G_4(q)$ are all nonnegative.
After a simple calculation, we get
\begin{align}\label{G3}
G_3(q)&=q^2\prod_{i=1}^{\infty}
\left(1+\frac{q^{2i}}{1-q^i}\right)\nonumber\\
&=q^2(1+q^2+q^3+2q^4+q^5+4q^6+2q^7
+6q^8+5q^9+9q^{10}+7q^{11}+\sum_{n=12}^{\infty}b_nq^n),\nonumber\\
&=q^2+q^4+q^5+2q^6+q^7+4q^8+2q^9+6q^{10}+5q^{11}+9q^{12}+7q^{13}
+\sum_{n=12}^{\infty}b_nq^{n+2},
\end{align}
where $b_n$ are nonnegative for $n\geq 12$.
Using \eqref{G1} and \eqref{G3}, we deduce that
\begin{align}\label{G1G3}
&G_1(q)+G_3(q)\nonumber\\
&\hskip 0.5cm=-q+q^2+q^4-q^5+2q^6-q^7+4q^8+q^9+6q^{10}+3q^{11}+9q^{12}+6q^{13}
+\sum_{n=12}^{\infty}b_nq^{n+2}.
\end{align}
Applying \eqref{G1G3} into \eqref{mainf4}, we have
\begin{align}\label{mainf5}
(1-q+q^4)\prod_{i=1}^{\infty}\left(1+\frac{q^{2i}}{1-q^i}\right)^2
&=G_1(q)+G_2(q)+G_3(q)+G_4(q)\nonumber\\[6pt]
&=-q+q^2+q^4-q^5+2q^6-q^7+4q^8+q^9+6q^{10}+3q^{11}\nonumber\\[5pt]
&\hskip 0.5cm+9q^{12}+6q^{13}
+\sum_{n=12}^{\infty}b_nq^{n+2}+G_2(q)+G_4(q).
\end{align}
Since the coefficients of $q$, $q^5$ and $q^7$ in $G_2(q)+G_4(q)$ are all zero,  the coefficients of $\eqref{mainpo}$ are nonnegative
except for $n=1,5,7$. This completes the proof.\qed

With the aid of Theorem \ref{thmuni} and  \eqref{pdsym}, we
arrive at the following corollary.
\begin{coro}
For $n\not= 1,5,7$, we have
\[
M_{bd}(n,n)\leq\cdots\leq M_{bd}(-1,n)\leq M_{bd}(0,n)
\geq M_{bd}(1,n)\geq\cdots \geq M_{bd}(n,n).
\]
\end{coro}
This means the sequence $\{M_{bd}(m,n)\}_{|m|\leq n}$
is unimodal for $n\not= 1,5,7$.

\noindent{\bf Acknowledgments.} The first author was supported by the Scientific Research Foundation of Nanjing Institute of Technology
(No. YKJ201627). The second author was supported by the Natural Science Foundation of Jiangsu Province of China (No. BK20160855) and the National
Natural Science Foundation of China (No. 11801139).

\end{document}